\documentclass[english,11pt]{amsart}
\usepackage{a4wide}
\usepackage{amsmath,amssymb,graphicx}
\usepackage{mathrsfs}
\usepackage{mathtools}
\usepackage{float}
\usepackage{dsfont}
\usepackage[all,cmtip]{xy}

\usepackage{tikz}
\usetikzlibrary{intersections,calc,arrows.meta, matrix}
\usetikzlibrary{arrows}
\usepackage{tikz}
   \usetikzlibrary{positioning}
\usepackage{bbding,wasysym}

\graphicspath{{figures/}}
\usepackage{hyperref}
\hypersetup{
	colorlinks=true,
	linkcolor=blue,    
}

\usepackage{frcursive}
\usepackage[utf8]{inputenc} 
\usepackage{amscd}

\def\di{\displaystyle}

\newtheorem{theorem}{Theorem}
\newtheorem{definition}{Definition}
\newtheorem{lemma}{Lemma}

\newtheorem{remark}{Remark}

\newcommand{\N}{\mathbb{N}}

\newcommand{\R}{\mathbb{R}}

\newcommand{\T}{\mathbb{T}}

\newcommand{\ack}{\section*{Acknowledgements}}
\begin{document}
\title{Mid-point embedding of Hamiltonian systems and variational integrators}
\author{Jacky CRESSON$^1$ and Rouba SAFI$^{1,2}$}
\begin{abstract}
Following the discrete embedding formalism, we give a new derivation of the mid-point variational integrators as developed by J.M. Wendlandt and J.E. Marsden by defining an adapted order two discrete differential and integral calculus. This allows us to obtain a clearer correspondence between the discrete and continuous case. We also discuss the corresponding definition of a discrete Hamiltonian system. A complete comparaison with the results of J.M. Wendlandt and J.E. Marsden is provided. 
\end{abstract}

\maketitle

$^1$ Laboratoire de mathématiques et leurs applications, UMR CNRS 5142, Université de Pau et des Pays de l'Adour-E2S.

$^2$ LaMA, Laboratoire de math\'ematiques et applications (Tripoli, Liban), Universit\'e libanaise. 

\tableofcontents

\section{Introduction}

In this article, we focus on the {\bf numerical integration of Hamiltonian systems}. This subject has leaded to a vast literature. The main point is that in order to obtain accurate and dynamically consistent numerical integrators one needs to preserve some classical properties of Hamiltonian systems. Some article are focusing on the preservation of energy, some other on the symplectic character of the flow or the variational structure of the equation. The preservation of these properties has lead to the field of {\bf Geometric numerical integration} \cite{hairer}.\\

In this paper, we focus on {\bf variational integrators} as they were defined by J.E. Marsden and M. West in \cite{marsden} and in particular to one of them called the mid-point variational integrators as it was discussed in \cite{wendlandt}. The main idea of variational integrators is to derive a discretisation of equations possessing a variational structure, i.e. whose solutions are corresponding to critical points of a functional, by a discretisation of the functional and a characterisation of the discrete critical points of this functional. As already said, this has already been done and the accuracy of resulting numerical integrators is well studied. Nevertheless, we believe that the formulation and the derivation of these integrators can be made more transparent following the formalism of {\bf discrete embedding} as exposed in \cite{cr1,cr2,cr3} by introducing a {\bf discrete differential and integral calculus}. Moreover, we thinks that this approach can be also interesting to role of each ingredients (quadrature formula, approximations) entering in the derivation of these integrators.\\

As an example, if one consider a classical Lagrangian systems whose functional is given by 
\begin{equation}
\label{lagfunc}
    \mathscr{L} (q)=\di\int_a^b L(q(s), \dot{q} (s) )\, ds ,
\end{equation}
where $L:\R^d \times \R^d \mapsto \R$ is denoted $L(q,v)$, the associated Euler-Lagrange equation is
\begin{equation}
\label{el}
\di\frac{d}{dt} \left (    \di\frac{\partial L}{\partial v} (q(t), \dot{q} (t) \right )
=\di\frac{\partial L}{\partial q} (q(t) ,\dot{q} (t) ),
\end{equation}
where $q:\R \mapsto \R^d$. The Marsden-West approach to these integrators in the most simple case consists in replacing the functional \ref{lagfunc} by the discrete analogue 
\begin{equation}
\mathscr{L} (q_0 ,\dots ,q_n ) =\di\sum_{i=0}^{n-1} \mathbb{L} (q_{i+1} ,q_i ) ,
\end{equation}
where 
\begin{equation}
    \mathbb{L} (q_{i+1} ,q_i )=h L\left ( q_i ,\di\frac{q_{i+1} -q_i}{h} \right ) ,
\end{equation}
for a small time increment $h$ using a first order quadrature formula for the integral and a first quadrature formula for the derivative. 

Denoting by $\mathbb{L}(y,x)$ the variables entering in the function $\mathbb{L}$, the resulting {\bf discrete Euler-Lagrange equation} is written as
\begin{equation}
\label{marsel}
\di\frac{\partial \mathbb{L}}{\partial x} (q_{i+1} ,q_i ) +\di\frac{\partial \mathbb{L}}{\partial y} (q_i ,q_{i-1} ) =0,\ i=1,\dots ,n-1 .     
\end{equation}
It is difficult to see directly how the discretisation acts on the initial Euler-Lagrange equation. The reason is that discussing the construction of a variational integrator using the point of view of the theory of approximation wich is dominant in numerical analysis, we forget about the dynamical nature of the Euler-Lagrange equation, meaning that such a Lagrangian functional is intimately related to the fact that position and speed are optimizing some quantities. \\

In order to recover the classical interpretation, one is lead to develop discrete analogue of the differential and integral calculus and in a second step a discrete analogue of the calculus of variations. This has been done up to now only in the context of a theory of order one, meaning that the discrete differential and integral calculus provide order one approximation to the continuous analogues. A complete presentation of the order one theory and a comparison with the classical work of J.E. Marsden and M. West in \cite{marsden} is provided in \cite{crouba}.\\

Let us explain briefly how it works. We first introduce two discrete operators acting on discrete function $f\in C(\T ,\R^d )$, where $\T =\{ t_i \}_{i=0,\dots ,n }$, $t_{i+1} -t_i =h$, $i=0,\dots ,n-1$ and $h>0$ is a small real increment.\\

The {\bf discrete right (resp. left) first-order derivative} denoted by $\Delta_+$ (resp. $\Delta_-$) is defined as 
\begin{equation}
\Delta_+ [f] (t_i )=\di\frac{f(t_{i+1} )-f(t_i )}{h} \ \ \ 
\left ( 
\mbox{\rm resp.}
\ \Delta_- [f] (t_i )=\di\frac{f(t_{i} )-f(t_{i-1} )}{h}
\right  ) 
,\ \ t_i \in \T^+ \ (\mbox{\rm resp.}
\T^- ) ,
\end{equation}
where $\T^+ =\T \setminus \{t_n \}$ (resp. $\T^- =\T \setminus \{ t_0 \}$). \\

The {\bf discrete integral} denoted by $\di\int_{t_i}^{t_j} f(s)\Delta_+ s $, $t_i ,\ t_j \in \T$, $j > i$, is defined as 
\begin{equation}
\di\int_{t_i}^{t_j} f(s)\Delta_+ s =\di\sum_{k=i}^{j-1} f(t_k ) h ,
\end{equation}
and $\di\int_{t_i}^{t_i} f(s)\Delta_+ s =0$ and $\di\int_{t_j}^{t_i} f(s)\Delta_+ s =-\di\int_{t_i}^{t_j} f(s)\Delta_+ s$.\\

These two operators mimics the classical differential and integral calculus. Precisely, we have again a {\bf fundamental theorem of the discrete differential calculus} given by 
\begin{equation}
    \di\int_a^b \Delta_+ {f} (s) \Delta_+ s = f(b)-f(a),\ \ \ \di\Delta_+ \left [ 
    \di\int_a^t f(s) \Delta_+ s 
    \right ] 
    = f(t)\ \forall \ t\in \T .
\end{equation}
A fundamental tools in the classical calculus of variations is the integration by parts formula. In the discrete case, it looks like
\begin{equation}
    \di\int_a^b f(s) \Delta_+ [g] (s) \Delta_+ s =-\di\int_a^b \Delta_- [f] (s) g(s) \Delta_+ s + f(b)g(b)-f(a)g(a) ,
\end{equation}
which put in evidence the duality between the left and right discrete derivative. \\

Having the discrete differential and integral calculus, the construction of a variational integrator is obtained following the {\bf discrete embedding formalism} developed in \cite{cr1,cr2,cr3}. We first define the {\bf discrete Lagrangian functional} denoted by $\mathscr{L}_h$ over $C(\T ,\R^d )$ by 
\begin{equation}
\mathscr{L}_h (q)=\di\int_a^b L (q(s) ,\Delta_+ [q] (s) ) \, \Delta_+ s , 
\end{equation}
for all $q\in C (T,\R^d )$.\\

A {\bf discrete calculus of variations} is developed as usual taking variations in the functional space 
\begin{equation}
    \mathscr{V} =\{ 
    f\in C(\T ,\R^d )\,, f(a)=f(b)=0 \} .
\end{equation}
The Frechet derivative $D\mathscr{L}_h (q)$ of $\mathscr{L}_h$ at point $q\in C(\T ,\R^d)$ in the direction $w\in \mathscr{V}$ is then given by 
\begin{equation}
    D\mathscr{L}_h (q) (w)=\di\int_a^b \left [ 
    \Delta_+ [w] (s) \di\frac{\partial L}{\partial v} (q(s),\Delta_+ [q] (s)) +w(s) \di\frac{\partial L}{\partial q} (q(s),\Delta_+ [q] (s))
\right ] 
\, \Delta_+ s
\end{equation}
which can be directly rewritten using the discrete integration by part formula and the boundaries conditions given in $\mathscr{V}$ as 
\begin{equation}
    D\mathscr{L}_h (q) (w)=\di\int_a^b -\left [ 
    \Delta_- \left [  \di\frac{\partial L}{\partial v} (q(s),\Delta_+ [q] (s)) \right ] +\di\frac{\partial L}{\partial q} (q(s),\Delta_+ [q] (s))
\right ] w(s)
\, \Delta_+ s
\end{equation}
The critical points of $\mathscr{L}_h$ satisfy $D\mathscr{L}_h (q) (w)=0$ for all $w\in \mathscr{V}$. As a consequence, we obtain the {\bf discrete Euler-Lagrange equation} as 
\begin{equation}
    -\Delta_- \left [  \di\frac{\partial L}{\partial v} (q(s),\Delta_+ [q] (s))
\right ] +\di\frac{\partial L}{\partial q} (q(s),\Delta_+ [q] (s)) =0 ,\ s\in \T^{\pm},
\end{equation}
where $\T^{\pm} =\T^+ \cap \T^-$. \\

This equation coincides with the Marsden's form given in \eqref{marsel} but we can notice at least two differences: \\

\begin{itemize}
\item The algebraic structure of the classical Euler-Lagrange equation is preserved. Moreover, one can see how the different discrete versions of the derivative and integral come into play in the derivation of the discrete version of the Euler-Lagrange equation.\\

\item The natural dichotomy between position and speed is also preserved in our formulation as the Lagrangian acts on a position $q$ and a discrete speed $\Delta_+ [q]$. This property will play a fundamental role in our derivation of a discrete version of Hamiltonian systems. In the contrary, Marsden's phase space is obtained by doubling the discrete configuration space then breaking the classical geometric approach.
\end{itemize}
\vskip 2mm
The benefit to develop a convenient discrete differential and integral calculus to present variational integrators is then clear but is restricted up to now to variational integrators of order one. \\

In this paper, {\bf we develop an order two discrete differential and integral calculus} which combined with discrete embedding formalism lead to a new formulation of the {\bf mid-point variational integrator} first defined by J.E. Marsden and J.M. Wendlandt in \cite{wendlandt} as a first step to the definition of a high-order differential and integral formalism.\\

A second problem appears when dealing specifically with the definition of what can be called a {\bf discrete Hamiltonian system}. Formally, we must answer the two following questions:\\

\begin{itemize}
    \item What is the discrete analogue of the phase space for Hamiltonian systems ? 
    
    \item What is the discrete definition of the Legendre transform ?
\end{itemize}

Indeed, as longs as the map $v\rightarrow \di\frac{\partial L}{\partial v} (q,v)$ is invertible for all $q$ in the configuration space, one can defined the Legendre transform denoted by $F$. This mapping acts on the tangent space of the configuration space and goes in the cotangent space.\\

In the Marsden's approach, the tangent space is replaced by a doubling of the configuration space $M\times M$ due to the fact that the equation is seen as a mapping depending on $(q_{i+1} ,q_i )$. In the discrete embedding case, we recover the classical approach as we preserve the notion of a derivative.\\

Second, the discrete Euler-Lagrange equation leads to the definition of a discrete momentum $p$ defined by 
\begin{equation}
    p(t)=\di\frac{\partial L}{\partial v} (q(t),\Delta_+ [q](t) ) .
\end{equation}
As a consequence, the discrete Hamiltonian system associated to the Lagrangian $L$ is given by the finite differences system
\begin{equation}
\label{ham2}
\left .
\begin{array}{llll}
    \di\Delta_- [p] & = & -\di\frac{\partial L}{\partial q}(q,p) , \ \ &  t\in \T^{\pm}\\
    \di\Delta_+ [q] & = & \mathbb{F} (q,p) & t\in \T^{\pm} ,
\end{array}
\right .
\end{equation}
where $\mathbb{F}$ is the inverse function of $\partial L/\partial v$. Introducing, the classical Hamiltonian function
\begin{equation}
    H(q,p)=p \mathbb{F} (q,p)-L(q,\mathbb{F} (q,p) )  ,
\end{equation}
corresponding to the Hamiltonian function in the continuous case, one can proved that the finite differences system \eqref{ham2} can be written as (see \cite{cp2}):
\begin{equation}
    \label{hameq}
\left .
\begin{array}{llll}
    \di\Delta_- [p] & = & -\di\frac{\partial H}{\partial q} (q,p) , \ \ &  t\in \T^{\pm}\\
    \di\Delta_+ [q] & = &  \di\frac{\partial H}{\partial p} (q,p) & t\in \T^{\pm} ,
\end{array}
\right .
\end{equation}
which preserves the classical form of continuous Hamiltonian systems. A finite differences system of the form \eqref{hameq} will be called a {\bf discrete (right\footnote{The word "right" corresponds to the fact that we have used the right discrete derivative in the definition of our functional associated with the corresponding discrete integral calculus.}) Hamiltonian system}.\\ 

It must be noted that this definition does not appears in this form in the work of Marsden's school. For example, in (\cite{hairer},VI.6.2) a discrete Hamiltonian system is defined from the discrete Euler-Lagrange equation \eqref{marsel} by choosing
\begin{equation}
    p_i =-\di\frac{\partial \mathbb{L}}{\partial y} (q_{i+1} ,q_i )
\end{equation}
which induces that (see \cite{hairer},p.207):
\begin{equation}
    p_{i+1} =\di\frac{\partial \mathbb{L}}{\partial x} (q_{i+1} ,q_i ) . 
\end{equation}
The resulting discrete Hamiltonian system is then described by the map $\Phi: (q_i ,p_i)\mapsto (q_{i+1}, p_{i+1} )$ (see \cite{lalwest}). As we have 
\begin{equation}
    \di\frac{\partial \mathbb{L}}{\partial x} (q_{i+1} ,q_i ) =\di\frac{1}{h}\di\frac{\partial L}{\partial v} (q_i ,\Delta_+ [q]_i ) ,
\end{equation}
we denote by $\mathbb{F}$ the inverse mapping defined for all $q_i$ and $p_{i+1}$ by 
\begin{equation}
    q_{i+1} = \mathbb{F} (q_i ,p_{i+1}) .
\end{equation}
The mapping $\Phi$ is symplectic (see \cite{hairer}, Theorem 6.1 p.207). However, the structure of the finite differences system destroys the specific Hamiltonian form of the continuous case due to the fact that choosing as coordinates $(q_{i+1}, q_i)$ to describe the functional instead of $(q_i ,v_i =\Delta_+ [q]_i )$  distorts the Legendre map $F$. 

\begin{remark}
It must be noted that our definition and construction coincides with the non-shifted definition of Hamiltonian systems on time scales as introduced by F. Pierret in \cite{p} after the work of C.D. Ahlbrandt in the discrete case \cite{ahlbrandt} and C.D. Ahlbrandt, M. Bohner and J. Ridenhour \cite{ridenhour} in the shifted case. 
\end{remark}

In this paper, {\bf we extend the previous construction and definition of order one discrete Hamiltonian systems to the case of order two using the order two discrete differential and integral calculus}.\\

The plan of the manuscript is as follows: \\

In Section \ref{discretediff}, we develop an order two discrete differential and integral calculus which coincide in the integral case with the mid-point quadrature formula. Different time-scales need to be introduced in order to take care of the quadrature formula and will induces several technical difficulties in particular for what concerns the integration by parts formula. Section \ref{midpointlag} deals with the definition of the mid-point Lagrangian functional and the associated calculus of variations. Having the mid-point Euler-Lagrange equation we can look for mid-point Hamiltonian systems. This is done in Section \ref{midpointham} where mid-point Legendre transform and momentum are defined. We prove that the solutions of the mid-point discrete Hamiltonian system correspond to critical points of the mid-point embedding of the classical Hamiltonian functional.Finally, Section \ref{conclusion} discuss several extension of the present work.

\section{Discrete mid-point Differential and integral calculus}
\label{discretediff}

\subsection{Definitions of different time scales}

Let us first set the definition of our different discrete time scales on $[a,b]$.

\begin{definition}\label{time scales}
Let $I=[a,b] \subset \R $,\, $N\in \N^*$ and let $h=(b-a)/N$, we define the following time scales
\begin{enumerate}
    \item[-]$\T=\{t_i=a+ih,\ i=0,1,...,N\}$.
    \item[-]$\T^+=\T \setminus \{b\},\;\;\T^-=\T \setminus \{a\} $. 
    \item[-]$\T_{\frac{1}{2}}=\{t_{i+\frac{1}{2}}=\frac{1}{2}(t_{i+1}+t_i),\;i=0,\dots,N-1\}$.
    \item[-]$\T_{\circ} =\T\cup\T_{\frac{1}{2}}$.
\end{enumerate}
\end{definition}

\begin{definition}
Let $\T$ a time scale defined on $[a,b]$. $\pi$ is a projection map on $\T_{\frac{1}{2}}$ defined by
\begin{align*}
    \pi \colon &\T^+ \to \T_{\frac{1}{2}}\\
    &t_i \mapsto \pi(t_i)=t_{i+\frac{1}{2}}=\frac{1}{2}\bigg[t_{i+1}+t_i\bigg].
\end{align*}
\end{definition}

\begin{definition}
let $\Tilde{\T}$  a discrete time scale on $[a,b]$, with a step $\Tilde{h}$, we donate by $\sigma_{\Tilde{\T}}$ and $\rho_{\Tilde{\T}}$ two maps defined by 
\begin{align*}
\sigma_{\Tilde{\T}} \colon \Tilde{\T}^+  &\to \Tilde{\T}^- \\
t &\mapsto \sigma_{\Tilde{\T}}(t)=t+\Tilde{h},
\end{align*}
and 
\begin{align*}
\rho_{\Tilde{\T}} \colon \Tilde{\T}^-  &\to \Tilde{\T}^+ \\
t &\mapsto \rho_{\Tilde{\T}}(t)=t-\Tilde{h}.
\end{align*}
\end{definition}

In the following, we use the simplified notations:
\begin{equation}
\sigma=\sigma_{\T}\;(resp.\,\rho=\rho_{\T}),\ \sigma_{\frac{1}{2}}=\sigma_{\T_\frac{1}{2}}\;(resp.\,\\ \rho_{\frac{1}{2}}=\rho_{\T_\frac{1}{2}}),\ 
 \sigma_\circ=\sigma_{\T_\circ}\;(resp.\,\rho_\circ=\rho_{\T_\circ}).
\end{equation}

\begin{figure}[H]
    \centering
    \includegraphics[scale=1.5]{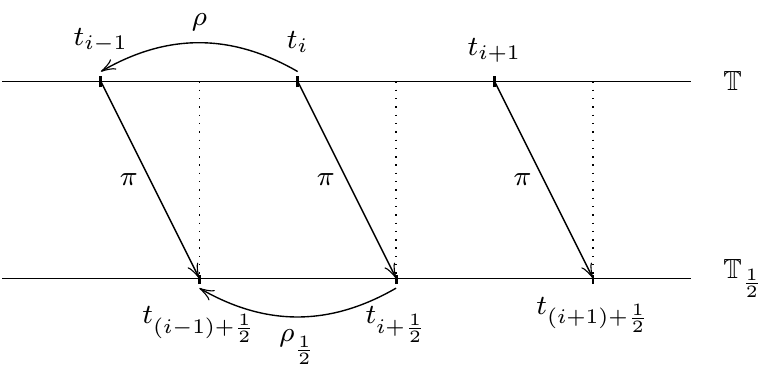}
    \caption{Connection between the different time scales}
\end{figure}

\subsection{Different functional spaces}

In the following, our basic objects are functions in $C(\T ,\R^d )$. However, the construction of the mid-point embedding uses functions over $C(\T_0 ,\R^d )$. These two functional spaces are connected via the following {\bf extension mapping}:

\begin{definition}
\label{extension} 
For all $f\in C(\T,\R^d)$ we define $f_\circ \in  C(\T_\circ,\R^d)$ and we call extension of $f$ on $\T_\circ$ the function defined by
\begin{equation}
f_\circ(t)=\;
\left\lbrace
\begin{array}{ll}
f(t_i),& t=t_i.\\\\
\dfrac{f(t_i)+f(t_{i+1})}{2},& t=t_{i+\frac{1}{2}}.
\end{array}
\right.
\end{equation}
\end{definition}

Another way to see the extension mapping is to introduce the {\bf interpolation map} denoted by $\iota_1$ over the set of piecewise continuous linear functions denoted by $\mathbb{P}^1$ of functions in $C(\T ,\R^d )$. Then, the mid-point extension $f_{\circ}$ of a given function $f$ is the projection over $\T_{\circ}$ of $\iota_1 [f]$. \\

In the following, we use different operators. 

\begin{definition}
    Let $f\in C(\T_{\frac{1}{2}} ,\R^d )$. We denote by $[f]_{\frac{1}{2},-}$ the function defined on $\T_{\frac{1}{2}}^-$ by
\begin{equation}
    [f]_{\frac{1}{2},-} (t)=\di\frac{1}{2} 
    \left ( 
    f(t)+f (\rho_{\frac{1}{2}} (t))
    \right ) ,
    \ \ \mbox{\rm for all}\ \ t\in \T_\frac{1}{2}^- .
\end{equation}
Equivalently, for $f\in C(\T_{\circ} ,\R^d )$, we denote by $[f]_{\circ}$ the function defined on $C (\T^{\pm} ,\R^d )$ by 
\begin{equation}
    [f]_{\circ } (t) =\di\frac{1}{2} \left (
    f(\sigma_{\circ} (t) )+f( \rho_{\circ} (t) )
    \right ) ,
    \ \ \mbox{\rm for all}\ \ t\in \T^{\pm} .
\end{equation}
\end{definition}

These manipulations will be useful when we will compute discrete integrals over $\T_{\frac{1}{2}}$ and interpreting them as discrete integrals over $\T$.\\

We introduce also the following notation:\\

For an arbitrary time scale $\tilde{\T}$ and a function $f$ of $C([a,b],\R^d )$ we denote by $\pi_{\tilde{\T}}$ the map from $C([a,b],\R^d )$ into $C(\tilde{\T},\R^d )$ obtained by taking the restriction of $f$ over $\tilde{\T}$.

\subsection{Discrete derivative and anti-derivative}

\subsubsection{Discrete derivatives}

Following the classical definition of derivatives on time scales as discussed in \cite{bohner}, discrete derivatives over an arbitrary discrete time scale are defined by:

\begin{definition}
Let $\tilde{\T}$ be an arbitrary time scale. We denote by $\Delta_{\tilde{\T},+}$ and $\Delta_{\tilde{\T},-}$ the operators defined $for\ all\ f\in (\tilde{\T},\R^d)$ by 
\begin{equation}
\Delta_{\tilde{\T},+}f(t)=\di\frac{f^{\sigma_{\tilde{\T}}}(t)-f(t)}{ \sigma_{\tilde{\T}}(t)-t} ,  
\end{equation}
and
\begin{equation}
\Delta_{\tilde{\T},-}f(t)=\di\frac{f(t)-f^{\rho_{\tilde{\T}}}(t)}{t-\rho_{\tilde{\T}}(t)} ,   
\end{equation}
with 
\begin{equation}
f^{\sigma_{\tilde{\T}}} := f\circ \sigma_{\tilde{\T}}\ \mbox{\rm and}\ \  f^{\rho_{\tilde{\T}}} :=f\circ \rho_{\tilde{\T}} .
\end{equation}
\end{definition}

It must be noted that the previous definition can be seen for $\tilde{\T} \in \{ \T , \T_{\frac{1}{2}} ,\T_{\circ} \}$ as follow:
\begin{equation}
\xymatrix{
    P^1_{\tilde{\T}} ([a,b],\R^d ) \ar[r]^{d^+ /dt}   &  P^{0,+}_{\tilde{\T}} ([a,b[,\R^d ) \ar[d]^{\pi_{\tilde{\T}}} \\
    C(\tilde{\T} ,\R^d )\ar[u]^{\iota_1}\ar[r]^{\Delta_{\tilde{\T},+}}   &  C(\tilde{\T}^+ ,\R^d )
}
,
\end{equation}
where $P^{0,+}_{\tilde{\T}}$ is the set of constant piecewise functions on intervals of the form $[\tilde{t}_i ,\tilde{t}_{i+1} [$ for $i=0,\dots ,\tilde{N}-1$ with $\tilde{t}_i,\tilde{t}_{i+1}$ in $\tilde{\T}$.\\

According to the time scale used, we simplify our notations as follows:\\

\begin{equation}
\left .
\begin{array}{ll}
\Delta_+ =\Delta_{\T,+} & \ (resp.\; \Delta_-=\Delta_{\T,-} ) .\\
\Delta_{\circ,+}=\Delta_{\T_\circ,+} & \ (resp.\; \Delta_{\circ,-}=\Delta_{\T_\circ,-} ).\\
\Delta_{\frac{1}{2},+}=\Delta_{\T_\frac{1}{2},+} & \ (resp.\; \Delta_{\frac{1}{2},-}=\Delta_{\T_{\frac{1}{2}},-} ).
\end{array}
\right .
\end{equation}

\subsubsection{Discrete anti-derivatives}

In the same way, one can define a discrete anti-derivative over an arbitrary discrete time scale (see \cite{bohner}): 

\begin{definition}
Let $\lambda \in [0,1[$. We denote by $t_{i,\lambda} =(1-\lambda )t_i +\lambda t_{i+1}$, $i=0,\dots ,N-1$. We denote by $\T_{\lambda}$ the set of $t_{i,\lambda}$, $i=0,\dots ,N-1$ and $\T_{\circ ,\lambda} =\T \cup \T_{\lambda}$. The $\lambda$-anti-derivative over $\T$ is defined for all function $f\in C(\T_{\circ ,\lambda}, \R^d )$ by 
\begin{equation}
\int_{t_i}^{t_{i+1}} f(t) \Delta_{\lambda , \T} t= f (t_{\lambda}) (t_{i+1} -t_i ).
\end{equation}
\end{definition}

It must be noted that despite the fact that we need the information about $\T_{\circ ,\lambda}$, the discrete anti-derivative is only defined on $\T$, meaning that we consider only integrals whose bounds of integration belong to $\T$.\\

We denote by $\iota_{0,\lambda ,\T,+}$ the mapping from $C(\T_{\lambda} ,\R^d )$ into $P^{0,+}_{\lambda,\T} ([a,b[,\R^d )$ defined by 
\begin{equation}
\iota_{0,\lambda ,\T, +} [f] (t)=
f(t_{i,\lambda} )\  \mbox{\rm for all} \  t\in [t_i ,t_{i+1} [,\ i=0,\dots ,N-1.
\end{equation}

The discrete anti-derivative is then obtained as follows:

\begin{equation}
\xymatrix{
    P^{0,+}_{\lambda,\T} ([a,b[,\R^d ) \ar[r]^{\di\int_a^t \dot\, ds}   &  P^{1,+}_{\T} ([a,b],\R^d ) \ar[d]^{\pi_{\T}} \\
    \ar[u]_{\iota_{0,\lambda, \T,+}} C(\T_{\lambda} ,\R^d )  \ar[r]^{\di\int_a^t \dot \, \Delta_{\lambda, \T} s} &   C(\T ,\R^d )
}
.
\end{equation}

We simplify our notations according to the time scales used. Precisely, we denote by 
\begin{equation}    
\displaystyle\int_{a}^{b} f(t) \Delta t=\displaystyle\int_{a}^{b} f(t) \Delta_{0,\T} t,\ 
    \displaystyle\int_{a}^{b} f(t) \Delta_{\frac{1}{2}} t=\displaystyle\int_{a}^{b} f(t) \Delta_{\frac{1}{2},\T} t .
\end{equation}

Using these notations, the classical {\bf mid-point quadrature formula} (see for example \cite{hairer}) for an integral of a function $f$ on $[a,b]$ over a discrete time scale $\T$ corresponds to the $1/2$-integral of the extension $f_{\circ}$ of $f_{/\T}$ over $\T_\frac{1}{2}$, i.e. 

\begin{equation}
\int_a^b f(s)\, ds \overset{\rm mid-point}{\simeq} \displaystyle\int_a^b f_{\circ} (t) \,\Delta_{\frac{1}{2},\T} t= \di\sum_{i=0}^{N-1} f_\circ (t_{i+\frac{1}{2}})(t_{i+1}-t_i)=\di\sum_{i=0}^{N-1} \di\frac{f(t_{i+1} )+f(t_i )}{2} (t_{i+1} -t_i ).
\end{equation}

\subsection{Proprieties of derivatives and integrals}

During the derivation of the mid-point Euler-Lagrange equation in Section \ref{midpointlag}, the computations mix objects coming from different time scales. As a consequence, we need to precise the connection between all these quantities.

\begin{lemma}
\label{relationprop}
For all $f\in C(\T,\R^d)$, we have 
\begin{equation}
\Delta_{\circ,+}[f_{\circ}](t_{i+\frac{1}{2}})=\Delta_+[f](t_i), \;\;\;\;\;for\, all\, i=0,...,N-1.
\end{equation}
\end{lemma}

Mixing of terms will naturally occurs in the discrete integration by part formula. Precisely, we have:

\begin{lemma}[Discrete integration by part formula]
\label{itegration by part}
let $f\in C(\T_\circ,\R^d)$ and $v\in C(\T,\R^d)$, we have
\begin{equation}
\displaystyle\int_{a}^{b} f(t)\Delta_{\circ,+}[v_{\circ} ](t) \, \Delta_{\frac{1}{2}} t=
-\displaystyle\int_{a}^{b} \Delta_{\frac{1}{2},-}[f]( \sigma_{\circ} (t) ) v(t)
\, \Delta t +f(t_{N-\frac{1}{2}})v(t_N)-f(t_{\frac{1}{2}})v(t_0)
.
\end{equation}
\end{lemma}

Two technical lemmas will be useful:\\

Let us denote by $C_0 (\T ,\R^d )$ the set \begin{equation}
     C_0(\T,\R^d)=\{v\in C(\T,\R^d);\, v(t_0)=v(t_N)=0\} .
\end{equation}

Then, we have:

\begin{lemma}\label{dubois1}
For all $f\in C(\T_\frac{1}{2},\R^d)$, $v\in C (\T,\R^d)$ we have 
\begin{equation}
\int_{a}^{b} f(t)v_\circ(t) \, \Delta_{\frac{1}{2}} t = 
    \int_{a}^{b} \left [ f \right ]_{\circ} (t) v(t) \, \Delta t + h \left( f(t_{N-\frac{1}{2}})v(t_N)+f(t_{\frac{1}{2}})v(t_0) \right),
\end{equation}
or equivalently 
\begin{equation}
\int_{a}^{b} f(t)\, v_{\circ} (t) \, \Delta_{\frac{1}{2}} t =  \int_{a}^{b} 
[ f ]_{\frac{1}{2},-} (\sigma_{\circ} (t)) v (t) \, \Delta t + h \left( f(t_{N-\frac{1}{2}})v(t_N)+f(t_{\frac{1}{2}})v(t_0) \right).
\end{equation}
Note that for $v \in C_0(\T, \R^d)$, the last term vanishes.
\end{lemma}
The proof is given in Section \ref{proofdubois1}.\\

As usual, an argument similar to the Dubois-Raymond lemma is needed. We remind the following result of the classical discrete calculus of variations:

\begin{lemma}[Discrete Dubois-Raymond lemma] Let $f\in C(\T ,\R^d )$ be a function such that 
\label{dubois}
\begin{equation}
    \int_{a}^{b} f(t) v(t) \Delta t=0,\;\; \mbox{\rm for all}\ v  \in C_0 (\T,\R^d),
\end{equation}
then $f(t)=0$ for all $t\in \T^{\pm}$.
\end{lemma}

The proof is given in Section \ref{proofdubois} for the convenience of the reader. 

\section{Discrete mid-point Lagrangian systems}
\label{midpointlag}

In this section, we use the mid-point discrete differential and integral calculus in order to associate to a given Lagrangian functional a discrete analogue. We follow the strategy of discrete embedding formalism as exposed in \cite{cr1,cr2,cr3}. We then develop the corresponding discrete calculus of variations and obtain a discrete mid-point Euler-Lagrange equation. Our result is compared with the work of J.M. Wendlandt and J.E. Marsden in \cite{wendlandt} about the same problem.

\subsection{Mid-point discrete Lagrangian functional}

A discrete functional is a mapping from $C (\T ,\R^d )$ in $\R$. A particular class of discrete functional are obtained by mid-point embedding of classical Lagrangian functionals:

\begin{definition}[Mid-point Lagrangian functionals]
Let $\T$ be a discrete time-scale on $[a,b]$. A discrete functional is called a Lagrangian functional if it exists a real valued function $L$ defined on $\R \times \R^d \times \R^d$ called the Lagrangian function such that 
\begin{equation}
\label{lagrange}
    \mathscr{L}_\T(q)=\int_{a}^{b} L(t,q_\circ(t),\Delta_{\circ,+}q_\circ(t)) \Delta_{\frac{1}{2}} t,
\end{equation}
for all $q\in C(\T,\R^d)$.
\end{definition}

In the following, we restrict our attention to Lagrangian functions which do not depend on time and we denote the variables by $(q,v)\in \R^d \times \R^d$.\\

The previous definition of a mid-point Lagrangian functional is fixed as long as the mid-point embedding is given. Formula \eqref{lagrange} gives the following explicit form for the mid-point Lagrangian functional:

\begin{align}
\label{explicitmidpointfunc}
\mathscr{L}_\T(q)&=\displaystyle\int_{a}^{b} L\bigg(q_\circ(t),\Delta_{\circ,+}q_\circ(t)\bigg) \Delta_{\frac{1}{2}}t\\
              &=\displaystyle\sum_{i=0}^{N-1} L\bigg(q_\circ(t_{i+\frac{1}{2}}),\Delta_{\circ,+}q_\circ(t_{i+\frac{1}{2}})\bigg)h\\
              &=\displaystyle\sum_{i=0}^{N-1} L\bigg(\dfrac{q(t_{i+1})+q(t_i)}{2},\dfrac{q(t_{i+1})-q(t_i)}{h}\bigg)h.
\end{align}

\subsection{Comparaison with the Wendlandt-Marsden discrete Lagrangian functional}

Our discrete Lagrangian functional \eqref{explicitmidpointfunc} coincides with the mid-point Lagrangian functional defined by Wendlandt and Marsden in \cite{wendlandt}. However, in their case, they do not introduce discrete analogues of the derivative and anti-derivative so that the complete analogy with the classical form of a Lagrangian functional is lost. Indeed, they introduce, starting from the quadrature formula \eqref{explicitmidpointfunc}, a new Lagrangian function $\mathbb{L}_h$ defined on $\R^d \times \R^d$ by 
\begin{equation}
\label{WML}
    \mathbb{L}_h (q_{i+1} ,q_i )= h\, L(\dfrac{q_{i+1} +q_i}{2},\dfrac{q_{i+1}-q_i}{h}) .
\end{equation}
It must be noted that this discrete Lagrangian corresponds in our setting to a discrete integral, namely for $q\in C(\T ,\R^d )$ we have
\begin{equation}
\mathbb{L}_h (q_{i+1} ,q_i ) = \di \int_{t_i}^{t_{i+1}} L(t,q_\circ(t),\Delta_{\circ,+}q_\circ(t)) \Delta_{\frac{1}{2}} t .
\end{equation}
As a consequence, the Wendlandt-Marsden discrete Lagrangian lead to the discrete function 
\begin{equation}
    \left .
    \begin{array}{lll}
    \T^+ & \rightarrow & \R ,\\
    t & \mapsto & \mathbb{L}_h (q^{\sigma} (t), q(t)) =\di \int_{t}^{\sigma (t)} L(t,q_\circ(t),\Delta_{\circ,+}q_\circ(t)) \Delta_{\frac{1}{2}} t .
    \end{array}
    \right .
\end{equation}
The Wendlandt-Marsden discrete Lagrangian functional is then given by 
\begin{equation}
    \mathbb{S} (q_0 ,\dots ,q_N )= \di\sum_{i=0}^{N-1} \mathbb{L}_h (q_{i+1},q_i ) .
\end{equation}
The term "functional" is not clear as $\mathbb{S}$ is a mapping from $\R^{d(N+1)}$ to $\R$. However, thanks to the one-to-one correspondence between the data of a $(q_0 ,\dots ,q_N ) \in \R^{d(N+1)}$ and $q\in C(\T ,\R^d )$ satisfying $q(t_0 )=q_0 ,\dots ,q(t_N ) =q_N$, we can introduce a discrete functional over $C(\T ,\R^d )$ denoted by $\mathscr{L}_{WM} (q)$ defined by 
\begin{equation}
\mathscr{L}_{WM} (q)=\mathbb{S} (q_0 ,\dots ,q_N ) .
\end{equation}
A direct computation shows that, for all $q\in C(\T,\R^d)$
\begin{equation}
\mathscr{L}_{WM} (q)=\mathscr{L}_{\T} (q) .
\end{equation}

In comparison, the classical diagram for the definition of a Lagrangian functional in the continuous case is 

$$\begin{tikzpicture}[node distance=5cm]
\begin{scope}
\node (A){
    $q$} ; 

\node [above=2cm of A](B){$(q, \dot{q})$};
\node [right=2.8cm of B](C){$L(q, \dot{q})$};
\node [right=2cm of A](D){$\mathscr{L}(q)=\di\int_{0}^{t}L(q, \dot{q})dt$};

\draw[->](A) --(B);
\draw[->](B) --node [above] {\scriptsize L}(C);
\draw[<-](D) --node [right] {\scriptsize $\di\int_{0}^{t}dt$ }(C);
\draw[->](A) --node [above] {\scriptsize $\mathscr{L}$ }(D);

\end{scope}
\end{tikzpicture}$$

Using our discrete differential and integral calculus, the discrete Lagrangian is given by 

$$\begin{tikzpicture}[node distance=5cm]
\begin{scope}
\node (A){
    $q$} ; 

\node [above=2cm of A](B){$(q_\circ, \Delta_{\circ,+}q_\circ)$};
\node [right=2.5cm of B](C){$L(q_\circ, \Delta_{\circ,+}q_\circ)$};
\node [right=2cm of A](D){$\mathscr{L}_\T(q)=\di\int_{0}^{t}L(q_\circ, \Delta_{\circ,+}q_\circ)\,\Delta_\frac{1}{2}t$};

\draw[->](A) --(B);
\draw[->](B) --node [above] {\scriptsize L}(C);
\draw[<-](D) --node [right] {\scriptsize $\di\int_{0}^{t}\Delta_\frac{1}{2}$ t}(C);
\draw[->](A) --node [above] {\scriptsize $\mathscr{L}_\T$ }(D);

\end{scope}
\end{tikzpicture}$$

We then recover a complete analogy between the continuous and the discrete case. \\

This correspondence is lost in the Wendlandt-Marsden case, precisely due to the fact that the understanding of the mapping $\mathbb{S}$ as a functional over discrete functions of $C(\T ,\R^d )$ is not used.\\

Another consequence, is the fact that there is no analogue of the mapping $q\mapsto (q,\dot{q} )$ in Wendlandt-Marsden contrary to the previous presentation, i.e. that classical Lagrangian function depends on two objects of different nature; namely position and speed. Here again, this is due to the fact that, as no discrete functions are used, no analogue of the derivative is described. This difference of point of views induces different conception of the phase space. Indeed, if $q$ has some unit $u$ then $\dot{q}$ has $u.t^{-1}$ as unit. However, in the Wendlandt-Marsden case, $\mathbb{S}$ is defined over quantities with the same unit, namely $(q_i ,q_{i+1} )$ breaking the signification of the mapping $q\mapsto (q,\dot{q})$. In our case, due to the mapping $q\mapsto (q_{\circ} ,\Delta_{\circ ,+} [q_{\circ} ] )$ we obtain quantities with unit $u$ and $u.t^{-1}$ as in the classical case.\\ 

The same phenomenon explain why the Wendlandt-Marsden discrete functional is not expressed explicitly as an integral over a discrete function.\\

The same remark will apply on the formulation of the mid-point Euler-Lagrange equation in the next Section.

\subsection{Discrete mid-point calculus of variations}

The discrete mid-point calculus of variations follows the usual construction of a discrete calculus of variations. We first specify the space of variations, i.e. the set of functions allowed during the deformation of the discrete Lagrangian functional.\\

We denote by $\mathscr{V}$ the set of variations defined by 

\begin{equation}
\mathscr{V} = \{ v \in C (\T,\R^d ),\, v(a)=v(b)=0 \} .
\end{equation}

We recover the usual set of variations for the order discrete calculus of variations.\\

The discrete Frechet derivative of $\mathscr{L}$ at point $q$ in the direction $v \in \mathscr{V}$ is given by :
\begin{equation}
D\mathscr{L}(q)(v)= \lim_{\epsilon \to 0} \frac{\mathscr{L}(q+\epsilon v)-\mathscr{L}(q)}{\epsilon}.
\end{equation} 

\begin{definition}
Let $\mathscr{L}$ a discrete Lagrangian functional. A critical point of $\mathscr{L}$ is a  discrete-time function $q \in  C(\T,\R^d )$, such that 
\begin{equation}
D\mathscr{L}(q)(v)=0,
\end{equation}
for all $v\in \mathscr{V}$, where $D\mathscr{L}(q)$ denotes the Frechet derivative of $\mathscr{L}$ at $q$.
\end{definition}

Our main result is the following Theorem:

\begin{theorem}[Discrete mid-point Euler-Lagrange equation]
Let $\mathscr{L}$ a discrete Lagrangian functional. The discrete Euler-Lagrange equation associated to $\mathscr{L}$ is given by

\begin{equation}
\label{ELPM}
    \bigg[ \frac{\partial L}{\partial q}( \star_{\circ} (t) ) \bigg] _ {\frac{1}{2},-}  
    = 
    \Delta_{\frac{1}{2},-} \bigg[ \frac{\partial L}{\partial v}(\star_{\circ} (t) ) \bigg] 
    ,\;\;\;\; \mbox{\rm for all}\ t\in\T_{\frac{1}{2}}^\pm,
\end{equation}
where we denote by $\star_{\circ} (t)$ the vector 
\begin{equation}
\star_{\circ} (t) =(q_\circ(t),\Delta_{\circ,+}[q_\circ](t) ) .
\end{equation}
\end{theorem}

The mid-point Euler-Lagrange equation induces a numerical scheme which enable us to determine $\star_{\circ} (t_{i+1+\frac{1}{2}})$ from the data of $\star_{\circ} (t_{i+\frac{1}{2}})$ for all $i=0,\dots,N-2$. Indeed, from the definition of $\Delta_{\circ,+}$, we have

\begin{equation}
\label{delta+32}
\Delta_{\circ,+}[q_0](t_{i+1+\frac{1}{2}}))=\frac{2}{h}q_0(t_{i+1+\frac{1}{2}})-\frac{2}{h}q_0(t_{i+\frac{1}{2}})-\Delta_{\circ,+}[q_0](t_{i+\frac{1}{2}}) ,
\end{equation}
the definition of $\Delta_{\frac{1}{2},-}$ and the mid-point Euler-Lagrange equation gives

\begin{equation}
\bigg[ \frac{\partial L}{\partial q}(
\star_{\circ} ) \bigg] _{\frac{1}{2},-} (t_{i+1+\frac{1}{2}} )
= 
\di\frac{1}{h} 
\left  (
    \frac{\partial L}{\partial v}(
    \star_{\circ} (t_{i+1+\frac{1}{2}}))
    -
    \frac{\partial L}{\partial v}(
    \star_{\circ} (t_{i+\frac{1}{2}}))
\right ),     
\end{equation}
as by definition, we have 
\begin{equation}
\left [
\frac{\partial L}{\partial q}(
\star_{\circ} ) \right ]
_{\frac{1}{2},-} (t_{i+1+\frac{1}{2}} ) = 
\di\frac{1}{2} 
\left ( 
\frac{\partial L}{\partial q}(
\star_{\circ} (t_{i+1+\frac{1}{2}} ))
+
\frac{\partial L}{\partial q}(
\star_{\circ} (t_{i+\frac{1}{2}} ))
\right ) .
\end{equation}
We obtain finally, replacing the quantity $\Delta_{\circ,+}[q_0](t_{i+1+\frac{1}{2}})$ by \eqref{delta+32} in the mid-point Euler-Lagrange equation for all $i=0,\dots ,N-2$, an implicit numerical scheme allowing us to determine $q_{\circ} (t_{i+1+\frac{1}{2}} )$. Using again formula \eqref{delta+32}, we then compute $\Delta_{\circ,+}[q_0](t_{i+1+\frac{1}{2}})$.\\

Of course, one is not interested in $q_{\circ}$ but in $q$. Using the quantities $q_{\circ} (t_{i+\frac{1}{2}})$ and $\Delta_{\circ,+}[q_{\circ}](t_{i+\frac{1}{2}})$, one can reconstruct $q$ as follows: 

\begin{equation}
\left ( 
\begin{array}{c}
q_{\circ} (t_{i+\frac{1}{2}}  )\\
\Delta_{\circ,+}[q_{\circ}](t_{i+\frac{1}{2}})
\end{array}
\right ) 
= 
\left ( 
\begin{array}{cc}
\frac{1}{2} & \frac{1}{2} \\
-\frac{1}{h} & \frac{1}{h}
\end{array}
\right ) 
\left ( 
\begin{array}{c}
q (t_i ) \\
q (t_{i+1} )
\end{array}
\right ) 
.
\end{equation}

As a consequence, the numerical scheme can be implemented as long as initial conditions $q_{\circ} (t_{\frac{1}{2}})$ and $\Delta_{\circ,+}[q_{\circ}](t_{\frac{1}{2}})$ are given or equivalently fixing the values of $q_0$ and $q_1$.

\subsection{Comparison with the Wendlandt-Marsden Euler-Lagrange equation}

This form of the mid-point Euler-Lagrange equation must be compared with the one obtained by Wendlandt and Marsden in \cite{wendlandt}.\\

For all $x,y\in \R^d$, let us denote by $\star_{x,y}$ the quantity
\begin{equation}
\star_{x,y} = 
\left ( 
\begin{array}{cc}
\di\frac{x+y}{2} , \di\frac{y-x}{h} 
\end{array}
\right ) 
\end{equation}
We have 
\begin{equation}
\label{connectMN}
\star_{q_i ,q_{i+1}} = (q_{\circ ,i+\frac{1}{2}} ,v_{\circ ,i+\frac{1}{2}} ),
\end{equation}
where $q_{\circ ,i+\frac{1}{2}} =q_{\circ} (t_{i+\frac{1}{2}} )$ and $v_{\circ ,i+\frac{1}{2}} =
\Delta_{\circ ,+} [q_{\circ} ] (t_{i+\frac{1}{2}} )$.\\

Using the discrete Lagrangian $\mathbb{L}_h$ defined in equation \eqref{WML}, Wendlandt and Marsden derive the following form for the discrete Euler-Lagrange equation: 
\begin{equation}
\label{WMEL}
    \partial_x \mathbb{L}_h (q_i ,q_{i+1} )+ \partial_y \mathbb{L}_h (q_{i-1} ,q_i )=0 .
\end{equation}
The previous form of the discrete Euler-Lagrange equation destroys the usual algebraic form of the classical Euler-Lagrange equation in the contrary to our presentation.\\

However, from a formal point of view, equation \eqref{WMEL} coincides with our mid-point discrete Euler-Lagrange equation \eqref{ELPM}.\\

Indeed, simple computations give 
\begin{equation}
    \begin{array}{lll}
   \partial_x \mathbb{L}_h (q_i ,q_{i+1}) & = & h\bigg[\di\frac{1}{2}\di\frac{\partial L}{\partial q} (\star_{q_i ,q_{i+1}} ) -\di\frac{1}{h} \frac{\partial L}{\partial v} (\star_{q_i ,q_{i+1}} )\bigg]. \\
    \partial_y \mathbb{L}_h (q_{i-1} ,q_i ) & = & h\bigg[\di\frac{1}{2}\di\frac{\partial L}{\partial q} (\star_{q_{i-1} ,q_i})+\di\frac{1}{h} \frac{\partial L}{\partial v} (\star_{q_{i-1} ,q_i}) \bigg].
    \end{array}
\end{equation}

Using \eqref{connectMN} and replacing $\partial_x \mathbb{L}_h (q_i ,q_{i+1})$ and $\partial_y \mathbb{L}_h (q_{i-1} ,q_i )$ by their expressions in \eqref{WMEL}, we obtain 

\begin{equation}
   \bigg( \frac{\partial L}{\partial q}(q_\circ(t_{i+\frac{1}{2}}),\Delta_{\circ,+}[q_\circ](t_{i+\frac{1}{2}}))\bigg)_{\frac{1}{2},-} 
   -
    \Delta_{-,\frac{1}{2}}\bigg(\frac{\partial L}{\partial v}(q_\circ(t_{i+\frac{1}{2}}),\Delta_{\circ,+}[q_\circ](t_{i+\frac{1}{2}}))\bigg)=0 .
\end{equation}

\subsection{Proof of the mid-point Euler-Lagrange equation}

The proof follows the continuous strategy. We first compute the Frechet derivative of  $\mathscr{L}(q)$:

\begin{theorem}
\label{Frechetder}
Let $v\in \mathscr{V}$. For all $q\in C(\T,\R^d)$, the Frechet derivative of $\mathscr{L}_\T(q)$ is given by
\begin{equation}
D\mathscr{L}_\T(q)(v)= \di\int_{a}^{b} 
\left (  
      \left [
      \frac{\partial{L}}{\partial{q}}(\star_{\circ})\right ]_{\frac{1}{2},-} (\sigma_{\circ}(t))-
    \Delta_{\frac{1}{2},-} 
    \left [
    \frac{\partial L}{\partial v} (\star_{\circ} )
    \right ] (\sigma_{\circ} (t)) 
\right ) v(t) \, \Delta t .
\end{equation}
\end{theorem}

The proof is given in Section \ref{proofrechet}.\\

A critical point of $\mathscr{L}_\T$ satisfies $D\mathscr{L}_\T(v)=0$ for all $v\in \mathscr{V}$. As a consequence, we have for all $v\in \mathscr{V}$
\begin{equation}
\di\int_{a}^{b} 
\left (  
      \left [
      \frac{\partial{L}}{\partial{q}}(\star_{\circ})\right ]_{\frac{1}{2},-} (\sigma_{\circ}(t))-
    \Delta_{\frac{1}{2},-} 
    \left [
    \frac{\partial L}{\partial v} (\star_{\circ} )
    \right ] (\sigma_{\circ} (t)) 
\right ) v(t) \Delta t =0.    
\end{equation}

Using the discrete Dubois-Raymond lemma \ref{dubois}, we conclude the proof. 

\subsection{Example: mid-point discretisation for Lagrangian from mechanics} 

We consider the classical class of Lagrangian from Mechanics which are of the form
\begin{equation}
L(q,v)=\frac{1}{2} v^2-V(q) ,
\end{equation}
where $(q,v)\in \R^d \times \R^d$ and the potential $V:\R^d \rightarrow \R$ is a $C^1$-function.\\

As $\di\frac{\partial L}{\partial v}(q,v)=v$, we have 
\begin{equation}
\frac{\partial L}{\partial v}(q_\circ,\Delta_{\circ,+}[q_\circ])=\Delta_{\circ,+}[q_\circ] .
\end{equation}

As a consequence, the mid-point Euler-Lagrange equation \eqref{ELPM} reads for $t=t_{i+\frac{1}{2}}\in \T_{\frac{1}{2}}^{\pm}$ as

\begin{equation}
         \frac{q(t_{i+1})-2q(t_i)+q(t_{i-1})}{h^2}=\frac{1}{2}\bigg[\frac{\partial
        L}{\partial q}\bigg(q_\circ(t_{i+\frac{1}{2}}),\Delta_{\circ,+}[q_\circ](t_{i+\frac{1}{2}})\bigg)+\frac{\partial
        L}{\partial q}(q_\circ(t_{(i-1)+\frac{1}{2}}),\Delta_{\circ,+}[q_\circ](t_{(i-1)+\frac{1}{2}}))\bigg] ,
\end{equation}
which is equivalent to 
\begin{equation}
        \frac{q(t_{i+1})-2q(t_i)+q(t_{i-1})}{h^2}=-\frac{1}{2}\bigg[\frac{\partial
        V}{\partial q}\bigg(\frac{q(t_{i+1})+q(t_i)}{2}\bigg)+\frac{\partial
        V}{\partial q}\bigg(\frac{q(t_{i})+q(t_{i-1})}{2}\bigg)\bigg]
\end{equation}

This last equation is the one obtained by Wendlandt and Marsden in \cite{wendlandt}.

\section{Discrete Mid-point Hamiltonian System}
\label{midpointham}

Having a definition of discrete mid-point Lagrangian systems, a natural question is to define the corresponding notion of discrete mid-point Hamiltonian systems. Following the discrete embedding strategy we define discrete mid-point momentum as the mid-point embedding of the classical continuous definition for momentum. However, and contrary to the order one case, this procedure is not trivial. This is due to the fact that the Legendre condition relates $p_{\circ}$ and $(q_0, \Delta_{\circ,+} [q_{\circ}] )$ and not directly $p$ as function of $C(\T^+ ,\R^d )$. As a consequence, we have a choice for the definition of $p$ to be made coherent with the mid-point embedding of the Legendre transform. This work leads us to a definition of discrete mid-point Hamiltonian systems which is very close to the continuous definition. We prove that it coincides with the definition proposed by Wendlandt and Marsden in \cite{wendlandt}.  

\subsection{Reminder about Hamiltonian systems}

\subsubsection{Definition and properties}

\begin{definition}
Let $d\in \N^\star $, $p\in C([a,b],\R^d)$ and $q\in C([a,b],\R^d)$. Let $H \colon \mathbb{R}^{2d} \to \mathbb{R}$. A $2d$-dimensional differential system of the form 
\begin{equation}
\label{HS}
\begin{pmatrix}
    dp\\
    dq
\end{pmatrix}
=J.\nabla H(p,q)dt,
\end{equation}

where

\begin{equation}
J=\begin{pmatrix}
0 & -Id\\
Id & 0
\end{pmatrix},\nabla H=\begin{pmatrix}
               \partial_pH\\
                \partial_qH
                \end{pmatrix}
,
\end{equation}

is called a Hamiltonian system with Hamiltonian $H$.
\end{definition}

An important property of Hamiltonian systems is that there solutions correspond to \emph{critical points} of a given functional, i.e. follows from a \emph{variational principle}.

\begin{theorem}[Variational Principle]
\label{principe variationelle}
The points $(p,q)\in C^1([a,b],\R^d) \times C^1([a,b],\R^d)$ satisfying Hamilton's equations are critical points of the functional
\begin{equation}
\mathscr{L}_H:\;{C^1([a,b],\R^d)\times C^1([a,b],\R^d)}\to\R
\end{equation}
defined by
\begin{equation}
\mathscr{L}_H(p,q) = 
\di\int_{a}^{b} \di L_H (p(t),q(t),\dot{p}(t),\dot{q}(t)) dt
\end{equation}
where $L_H:{\R^d \times \R^d \times \R^d \times \R^d}\to\R$ is the Lagrangian defined by
\begin{equation}
\label{actionfunc}
\di L_H(x,y,w,v)=x\cdot v-H(x,y).
\end{equation}
\end{theorem}

\subsubsection{Lagrangian versus Hamiltonian systems}	

As we have already a consistent discrete theory of Lagrangian system, we will use this derivation of Hamiltonian systems in order to define \emph{discrete Hamiltonian systems}. Doing so, we will see that all the objects are related by the discrete embedding procedure.\\
	
	Let $L$ be a Lagrangian and denote by \mbox{EL} the corresponding Euler-Lagrange equation given by 
	\begin{equation}
		\di \frac{d}{dt} 
		\left [ 
		\frac{\partial L}{\partial v} (q,\dot{q}) 
		\right ] 
		=\di\frac{\partial L}{\partial q} (q,\dot{q}) .
		\label{EL}
	\end{equation}
	We assume that the Lagrangian is \emph{admissible}, i.e. that the map $v\rightarrow \di\frac{\partial L}{\partial v}$ is invertible for all $(q,t) \in \R^n \times \R$. As a consequence, we can introduce the \emph{moment} variable
	\begin{equation}
		p=\frac{\partial L}{\partial v} (q,\dot{q}) ,
	\end{equation}
	in order to rewrite the Euler-Lagrange equation as a first order system of differential equations given by
	\begin{gather}
		\left .
		\begin{split}
			\dot{q} & =  g(p,q) ,\\
			\dot{p} & =  \di\frac{\partial L}{\partial q} (q, g(p,q)) ,
		\end{split}
		\right .
		\label{change}
	\end{gather}
	where $g$ is the inverse of $\di\frac{\partial L}{\partial v}$. This change of variable will be called \emph{Legendre transform} in the following. Introducing the Hamiltonian function 
	\begin{equation}
		\label{ham}
		H(p,q) =pg(p,q)- L(q,g(p,q)),
	\end{equation}
one proves that equation (\ref{change}) is Hamiltonian with respect to $H$

\begin{align*}
\frac{\partial H}{\partial p}(p,q)&=g(p,q)+p\frac{\partial g}{\partial p}(p,q)-\frac{\partial L}{\partial v}(q,\dot{q},t)\frac{\partial g}{\partial p}(p,q)=g(p,q)=\frac{dq}{dt}.\\
\frac{\partial H}{\partial q}(p,q)&=p\frac{\partial g}{\partial q}(p,q)-\frac{\partial L}{\partial q}(q,\dot{q})-\frac{\partial L}{\partial v}(q,\dot{q})\frac{\partial g}{\partial q}(p,q)=-\frac{\partial L}{\partial q}(q,\dot{q})=-\frac{dp}{dt}.
\end{align*}

\subsection{Toward discrete Hamiltonian systems}

In the classical discrete (order one) case, the natural definition of a discrete Hamiltonian system associated to a given discrete Euler-Lagrange equation is coherent, meaning that the critical point of the discrete embedding of the classical action functional corresponds to the writing of the discrete Euler-Lagrange equation using the Legendre transform. 

\subsubsection{Discrete momentum and discrete Legendre transform}

Let $\mathscr{L}$ a Lagrangian functional and let $\mathscr{L}_\T$ the discrete Lagrangian functional associated to $\mathscr{L}$ defined in \eqref{lagrange}.
Regarding to the discrete mid-point Euler-Lagrange equation and following the usual way to derive Hamiltonian system in the continuous case, it is natural to introduce the following definition of a {\bf discrete momentum}:

\begin{definition}[Discrete momentum constraint]
Let $L$ be a Lagrangian system. We call discrete momentum a function $p\in C(\T^+,\R)$ such that 
\begin{equation}
\label{momentdef}
p_\circ=\displaystyle\frac{\partial L}{\partial v}(q_\circ,\Delta_{\circ,+}[q_\circ]).
\end{equation}
\end{definition}

As a consequence, assuming that the function $L$ is {\bf admissible}, i.e. that for all $q\in \R^d$ the mapping $v\mapsto \di\frac{\partial L}{\partial v} (q,v)$ is invertible, and denoting by $g$ the inverse, we obtain 

\begin{equation}
   \Delta_{\circ,+}[q_\circ]=g(p_\circ,q_\circ) .
\end{equation}

A main question is to be able to construct a function $p$ satisfying condition \eqref{momentdef}. The specific form of this relation implies that we must have a relation of the form
\begin{equation}
\left \{
\begin{array}{lll}
p(t_{i+1}) & = & \partial_v L (\star_{t_{i+\frac{1}{2}}} ) +\di\frac{h}{2} w(t_{i+\frac{1}{2}} ) ,\\    
p(t_i) & = & \partial_v L (\star_{t_{i+\frac{1}{2}}} ) -\di\frac{h}{2} w(t_{i+\frac{1}{2}} ) 
,
\end{array}
\right .
\end{equation}
where $w$ is a function to be determined.\\

In order that the previous relations induce a coherent definition for the function $p$, we must have
\begin{equation}
p(t_i )=\partial_v L (\star_{t_{i-\frac{1}{2}}} ) +\di\frac{h}{2} w(t_{i-\frac{1}{2}} ) .
\end{equation}

The mid-point Euler-Lagrange equation can be used to precise a suitable function $w$. Indeed, we must have 

\begin{equation}
\partial_v L (\star_{t_{i+\frac{1}{2}}} ) - \partial_v L (\star_{t_{i-\frac{1}{2}}} )
= \di\frac{h}{2} 
\left [ 
\partial_q L (\star_{t_{i+\frac{1}{2}}} ) +
\partial_q L (\star_{t_{i-\frac{1}{2}}} ) 
\right ] .
\end{equation}
As by definition of $p$ we have 
\begin{equation}
\partial_v L (\star_{t_{i+\frac{1}{2}}} ) - \partial_v L (\star_{t_{i-\frac{1}{2}}} ) = \di\frac{h}{2} \left [ 
w(\star_{t_{i+\frac{1}{2}}})+w(\star_{t_{i-\frac{1}{2}}}) \right ] ,
\end{equation}
we deduce that a suitable choice for $w$ is
\begin{equation}
w (\star_{t_{i+\frac{1}{2}}})= \partial_q L (\star_{t_{i+\frac{1}{2}}}) .
\end{equation}

We then are leaded to the following definition of the discrete momentum:

\begin{definition}[Discrete momentum] 
\label{defdiscmom}
We call discrete momemtum associated to $L$ the discrete function $p\in C(\T^+ ,\R^d )$ defined by 
\begin{equation}
\label{mom1}
p(t_i )=\partial_v L (\star_{t_{i-\frac{1}{2}}}) +\di\frac{h}{2} \partial_q L (\star_{t_{i-\frac{1}{2}}} ) ,    
\end{equation}
for all $t_i \in \T^{\pm}$ and 
\begin{equation}
\label{mom2}
p(t_0 ) =\partial_v L (\star_{t_{\frac{1}{2}} }) -\di\frac{h}{2} \partial_q L (\star_{t_\frac{1}{2}} )
.
\end{equation}
\end{definition}

\subsubsection{Comparaison with Wendlandt and Marsden}

In \cite{wendlandt}, Wendlandt and Marsden take as a definition for the discrete momentum the quantities
\begin{equation}
p (t_i ) = -\partial_x \mathbb{L}_h (\star_{q_i,q_{i+1}} )\ \ \mbox{\rm and}\ \ 
p(t_{i+1} ) =\partial_y \mathbb{L}_h (\star_{q_{i-1},q_i} ) .
\end{equation}

Of course, the previous definition is not usual and far from the standard definition of the momentum in the continuous case. The choice of the minus sign is also not explained. The main remark is of course that in order that the previous equalities make sense, then we must be sure that taking as a definition for all $t_i \in \T^+$ the definition of $p(t_i )$ we must have 

\begin{equation}
p(t_{i+1} ) =- \partial_x \mathbb{L}_h (\star_{q_{i+1},q_{i+2}} ) 
,
\end{equation}
ensuring the coherence of the definition of $p$. This is of course the case using the fact that we are looking for discrete functions $q$ which are solutions of the mid-point Euler-Lagrange equation.\\

Expliciting the Lagrangian $\mathbb{L}_h$, we recover the formula given in definition \ref{defdiscmom}. \\

We can notice that $p_{\circ}$ can also be computed directly in the Wendlandt and Marsden case and coincide with our choice of a discrete function $p$ satisfying the discrete momentum constraint.

\subsection{Discrete Hamiltonian function and discrete Hamiltonian systems}

Following the usual strategy, we consider discrete Hamiltonian function associated to $L$:

\begin{definition}
The discrete mid-point Hamiltonian function associated to the Lagrangian $L$ is the mid-point embedding of the classical continuous Hamiltonian function, i.e. for all $(p_{\circ} ,q_{\circ} )$, we consider the discrete function 
\begin{equation}
\label{hamiltonian discret}
H(p_\circ,q_\circ) = -L(q_\circ,g(p_\circ,q_\circ)) +p_\circ g(p_\circ,q_\circ) .
\end{equation}
\end{definition}

Here again, as for the definition of the discrete mid-point momentum, the discrete analogue is obtained directly just taking the definition of the classical continuous function in the discrete framework. \\

Using this function, the discrete Euler-Lagrange equation \eqref{ELPM} can be rewritten as the following discrete system:

\begin{definition}
[Discrete mid-point Hamiltonian systems]
Let $L$ be an admissible Lagrangian, then the mid-point discrete Euler-Lagrange equation \eqref{ELPM} can be written as 
\begin{equation}
      (SH_{\circ})\;\;
      \left\lbrace
      \begin{array}{lllll}
         \Delta_{\frac{1}{2},-}[p_\circ]&(t)&=
         \bigg [ - \frac{\partial H}{\partial q}(p_\circ,q_\circ)\bigg]_{\frac{1}{2},-} (t), &&for\; all \;t\in \T_{\frac{1}{2}}^{+,-}.\\
         \\
          \Delta_{\circ,+}[q_\circ]&(t)&=\frac{\partial H}{\partial
          p}(p_\circ,q_\circ)(t), &&for\; all\; t\in \T_{\frac{1}{2}}^{+,-}.
    \end{array}
    \right.
\end{equation}
where  $(q,p) \in C(\T ,\R^d ) \times C(\T^+ ,\R^d )$, the discrete momentum $p$ satisfies \eqref{mom1} and \eqref{mom2} and $H$ is defined by \eqref{hamiltonian discret}.\\

A discrete system of the form $(SH_{\circ})$ is called a discrete mid-point Hamiltonian system.
\end{definition}

The implementation of the algorithm goes as follows. Choose an initial condition $(p_0 ,q_0 )$. Then, by assumption on the form of $p$ we have 
\begin{equation}
p_1 =p_0 +h \partial_q L \left ( 
\di\frac{q_0+q_1}{2} ,\di\frac{q_1 -q_0}{h} 
\right ) ,
\end{equation}
and by the Legendre relation
\begin{equation}
q_1=q_0 +h g
\left ( 
\di\frac{p_1 + p_0}{2} , \di\frac{q_0+q_1}{2} \right).
\end{equation}
As a consequence, knowing $(p_0 ,q_0)$ we can determine $(p_1 ,q_1 )$ and then $(q_{\circ} (t_{\frac{1}{2}}), p_{\circ} (t_{\frac{1}{2}}))$.\\

The discrete mid-point Hamiltonian system then determines the quantities $(q_{\circ} (t_{i+\frac{1}{2}} )$ and $p_{\circ} (t_{i+\frac{1}{2}} ) )$ recursively and as a consequence, the quantities $q_{i+1}$ and $p_{i+1}$ for $i\geq 1$.\\ 

It seems reasonable to take the previous system as a definition for a discrete Hamiltonian system. However, in order that this definition mimics the continuous case, it is suitable that the solutions of this discrete system are in correspondence with critical points of a suitable discrete action functional, the most natural choice being the discrete mid-point embedding of the classical continuous action functional. We discuss this problem in details in the next Section.

\subsection{A variational approach to discrete Hamiltonian systems}

As for the discrete mid-point Lagrangian functional, we define the discrete mid-point action functional using the discrete embedding strategy. Precisely, we have:

\begin{definition}[Discrete action functional]
The discrete action functional associated to the given discrete Hamiltonian system is defined by 
\begin{equation}
\label{L-H}
\mathscr{L}_{H,\T} (p,q) =\di\int_{a}^{b}
\left ( 
p_{\circ} \Delta_{\circ,+} [q_\circ]-H(p_\circ,q_\circ) \right ) 
\, \Delta_{\frac{1}{2}} t.
\end{equation}
\end{definition}

A natural question is wheather the solutions of the discrete mid-point Hamiltonian system coincides with the critical point of $\mathscr{L}_{H,\T}$.\\ 

\begin{definition}
A couple of functions $(q,p) \in C(\T ,\R^d ) \times C(\T^+ ,\R^d )$ is a critical point of the discrete action functional $\mathscr{L}_{H,\T}$ if and only if for all variations $(v,w)\in \mathscr{V}\times C(\T ,\R^d )$, we have 
\begin{equation}
D\mathscr{L}_{H,\T} (p,q)(v,w)=0 .
\end{equation}
\end{definition}

Note that there is no constraints on the variations associated to $p$. As a consequence, we can not apply directly the mid-point Euler-Lagrange equation which was derived for variations in $\mathscr{V}$. However, simple computations lead to:

\begin{theorem}
\label{criticalham}
The critical points of the discrete action functional \eqref{L-H} corresponds to the solutions of the mid-point discrete Hamiltonian system ($SH_{\circ}$). 
\end{theorem}

We then obtain a global coherent picture using the mid-point discrete embedding of the relation between discrete Lagrangian, discrete Hamiltonian, discrete Legendre transform under the mid-point discretisation.\\

The proof of Theorem \ref{criticalham} is given in Section \ref{proofham}.

\section{Conclusion and perspectives}
\label{conclusion}

Following the discrete embedding formalism, We have develop a discrete differential and integral calculus allowing us to include mid-point classical approach to the construction of variational integrators and the definition of mid-point Hamiltonian systems. This formalism allows us to obtain a direct understanding of the way the classical continuous Euler-Lagrange equation or Hamiltonian systems are extended under a discretization framework. In particular, it allows us to preserve the differential form as well as the integral form of all the continuous object (Euler-Lagrange equation, Lagrangian functional). For Hamiltonian systems, it gives a direct constraint on the definition of a discrete momentum and a form of mid-point Hamiltonian systems which preserves the classical continuous structure. These constructions offers a new understanding of the classical approach of Wendlandt and J.E. Marsden on the construction of mid-point variational integrators for which we provide a full comparaison with our results. \\

This work can be generalize in many directions:\\

\begin{itemize}
\item All the remarks made about the structure of the discrete Euler-Lagrange equation, the definition of the discrete phase-space and the definition of a discrete Hamiltonian system extend to the high order variational integrators version as derived for example first by M. Leok \cite{leok} (see also C.N. Campo et al. in \cite{campo}) using the approach of J.E. Marsden and M. West. The extension of the discrete differential and integral calculus in order to cover these cases is under progress \cite{cka}.

\item Mid-point quadrature are unovoidable in a first approach to Stratonovich sochastic integrals (see \cite{oksendal}). Moreover, stochastic Hamiltonian systems as defined by J-M. Bismut in \cite{bismut} (see also \cite{ortega}) are extension of the classical notion of Hamiltonian systems using the Stratonovich formalism. Stochastic variational integrators were already develop in this setting following the J.E. Marsden and M. West approach for example by N. Bou-Rabee et al. in \cite{rabee} or Wang et al. in \cite{wang}. A derivation of stochastic variational integrators for stochastic Hamiltonian systems in the framework of discrete mid-point differential and integral will be given in a future work.
\end{itemize}

\ack J. Cresson thanks the GDR CNRS no. 2043, Géométrie différentielle et Mécanique. R. Safi thanks the LASeR (Tripoli, Liban) for the grants "Bourse d'excellence". Both thank M. Jazar for support. 

\begin{appendix}

\section{Proof of the lemmas}

\subsection{Proof of Lemma \ref{relationprop}}

Let $f\in C(\T,\R^d)$ than 
\begin{align*}
\Delta_{\circ,+}[f_\circ](t_{i+\frac{1}{2}}) 
     &=\dfrac{f_\circ(t_{i+1})-f_\circ(t_{i+\frac{1}{2}})}{h/2}\\
     &=\dfrac{f(t_{i+1})-\frac{1}{2} [f(t_i)+f(t_{i+1})]}{h/2}\\
     &=\dfrac{f(t_{i+1})-f(t_i)}{h}=\Delta_+[f](t_i),
\end{align*}
for all $i=0,..., N-1$.

\subsection{Proof of the discrete integration by part formula}

For all $f\in C(\T_0,\R^d)$ and $v\in C(\T,\R^d)$ we have 
\begin{equation}
\begin{array}{ll}
       \displaystyle\int_{a}^{b} f(t)\Delta_{\circ,+}[v_\circ](t) \, \Delta_{\frac{1}{2}} t &=h\displaystyle\sum_{i=0}^{N-1} f(t_{i+\frac{1}{2}})\Delta_{\circ,+}[v_\circ](t_{i+\frac{1}{2}})\\
       
      &=h\displaystyle\sum_{i=0}^{N-1} f(t_{i+\frac{1}{2}}) \Bigg[\frac{v(t_{i+1})-v(t_i)}{h}\bigg]\\
      
      &=\displaystyle\sum_{i=1}^{N} f(t_{(i-1)+\frac{1}{2}})g(t_i)-\displaystyle\sum_{i=0}^{N-1} f(t_{i+\frac{1}{2}})g(t_i)\\
      
      &=h\displaystyle\sum_{i=1}^{N-1}
      \bigg[\frac{f(t_{(i-1)+\frac{1}{2}})-f(t_{i+\frac{1}{2}})}{h}\bigg]v(t_i)+f(t_{(N-1)+\frac{1}{2}})v(t_N)-f(t_{\frac{1}{2}})v(t_0)\\
      
      &=-\displaystyle\int_{a}^{b} \Delta_{\frac{1}{2},-}[f](\sigma_\circ(t))v(t) \,
      \Delta t +f(t_{N-\frac{1}{2}})v(t_N)-f(t_{\frac{1}{2}})v(t_0).
      \end{array}
\end{equation}

\subsection{Proof of Lemma \ref{dubois1}}
\label{proofdubois1}

Let $f\in C(\T_\frac{1}{2},\R^d)$ and $v\in C_\circ(\T,\R^d)$. By definition of $\Delta_\frac{1}{2}$ integral, we have
\begin{equation}
    \int_{a}^{b} f(t)v_\circ(t) \Delta_{\frac{1}{2}} t 
     =h \di\sum_{i=0}^{N-1} f(t_{i+\frac{1}{2}} ) v_{\circ} (t_{i+\frac{1}{2}} ) ,
\end{equation}
As $v_{\circ} (t_{i+\frac{1}{2}}) =(v(t_{i+1} ) +v(t_i ))/2$, we obtain, regrouping the terms
\begin{equation}
\left .
\begin{array}{ll}
\di\int_{a}^{b} f(t)v_\circ(t) \Delta_{\frac{1}{2}} t
    & =h\di\sum_{i=1}^{N-1} \left ( \di\frac{f(t_{i+\frac{1}{2}})+f(t_{i-\frac{1}{2}})}{2} \right ) v(t_i ) \\
    & +h \left ( f(t_{(N-1)+\frac{1}{2}})v(t_N ) + f(t_{\frac{1}{2}})v(t_0) \right )  .
\end{array}
\right .
\end{equation} 
By definition, we have for $i=1,\dots ,N-1$
\begin{equation}
\left [ f \right ]_{\circ} (t_i )= \di\frac{1}{2} \left [ f(t_{i+\frac{1}{2}})+f(t_{i-\frac{1}{2}}) \right ] ,
\end{equation}
so that the first sum can be written as a classical discrete integral
\begin{equation}
\int_{a}^{b} f(t)v_\circ(t) \Delta_{\frac{1}{2}} t = 
    \int_{a}^{b} \left [ f \right ]_{\circ} (t) v(t)\Delta t
+h \left ( f(t_{N-\frac{1}{2}})v(t_N ) + f(t_\frac{1}{2})v(t_0) \right )  .
\end{equation}
This concludes the proof.

\subsection{Proof of Lemma \ref{dubois}}
\label{proofdubois}

The proof follows the classical one. Indeed, we can choose $v\in C_0 (\T,\R^d)$ such that for all $t\in \T^\pm $ we have $v(t)=f(t)$ and $v(t_0 )=v(t_N )=0$ by assumption. Then, we have
\begin{equation}
\int_{a}^{b} f(t) v(t) \Delta t = \sum_{k=1}^{N-1} f(t_k ) v (t_k) h =h \sum_{k=1}^{N-1}f(t_k)^2 =0.
\end{equation} 
Therefore $f(t)=0$ for all $t\in \T^{\pm}$. 

\subsection{Proof of Theorem \ref{Frechetder}}
\label{proofrechet}

Let $v\in \mathscr{V}$. Denoting $\star_{\circ} =(q_\circ ,\Delta_{\circ,+}[q_\circ] )$.
\begin{equation}
\mathscr{L}_\T(q+\epsilon v)=
     \int_{a}^{b} L  \bigg(  (q+\epsilon v)_\circ(t),\Delta_{\circ,+}[(q+\epsilon v)_\circ](t)  \bigg) \, \Delta_{\frac{1}{2}}t .
\end{equation}
According to the linearity of the extension map given in Definition \ref{extension} and the linearity of the discrete derivative $\Delta_{\circ,+}$, we have
\begin{equation}
\mathscr{L}_\T(q+\epsilon v)=
    \int_{a}^{b} L\bigg(q_\circ(t)+\epsilon v_\circ(t),\Delta_{\circ,+}[q_\circ](t)+\epsilon \Delta_{\circ,+}[v_\circ](t)\bigg) \, \Delta_{\frac{1}{2}}t .
\end{equation}
Doing a Taylor expansion of $L$ around $\star_{\circ} (t)$, we obtain
\begin{equation}
D\mathscr{L}_\T(q)(v)=\di\int_{a}^{b} \frac{\partial{L}}{\partial{q}}(\star_\circ(t)) \, v_\circ(t)+ \frac{\partial{L}}{\partial{v}}(\star_{\circ} (t))\, \Delta_{\circ,+}[v_\circ](t)\, \Delta_{\frac{1}{2}}t .
\end{equation}
Using the discrete integration by part formula for $v\in \mathscr{V}$ and Lemma \ref{dubois1}, 
\begin{equation}
D\mathscr{L}_\T(q)(v)= \di\int_{a}^{b} 
\left (  
      \left [
      \frac{\partial{L}}{\partial{q}}(\star_{\circ})\right ]_{\frac{1}{2},-} (\sigma_{\circ}(t))-
    \Delta_{\frac{1}{2},-} 
    \left [
    \frac{\partial L}{\partial v} (\star_{\circ} )
    \right ] (\sigma_{\circ} (t)) 
\right ) v(t) \, \Delta t .
\end{equation}
This concludes the proof.

\subsection{Proof of theorem \ref{criticalham}}
\label{proofham}

By definition of $\mathscr{L}_{H,\T}$, we consider the discrete mid-point Lagrangian functional 

\begin{equation}
\mathscr{L}_{H,\T}(p,q)
= \di\int_{a}^{b} 
\mathbb{L}
\left (
p_\circ ,q_\circ ,\Delta_{\circ,+}[p_\circ],\Delta_{\circ,+}[q_\circ]
\right ) 
\Delta_{\frac{1}{2}} t
,
\end{equation}
where 
\begin{equation}
\mathbb{L} (p,q,w,v)=pv-H(p,q) .
\end{equation}

Computing the Frechet derivative of $\mathscr{L}_{H,\T}$ in the direction $(\tilde{p}, \tilde{q})$, we obtain
\begin{equation}
\left .
\begin{array}{lll}
D\mathscr{L}_{H,\T} (p, q) (\tilde{p},\tilde{q})
& = & \di\int_a^b 
\left ( 
\di\frac{\partial \mathbb{L}}{\partial q} (\star_{\circ} (t)) \tilde{q}_{\circ} (t) 
+ \di\frac{\partial \mathbb{L}}{\partial v} (\star_{\circ} (t)) \Delta_{\circ,+} [\tilde{q}_{\circ}] (t) 
\right ) 
\,\Delta_{\frac{1}{2}}t \\

& & + \di\int_a^b 
\left ( 
\di\frac{\partial \mathbb{L}}{\partial p} (\star_{\circ} (t)) \tilde{p}_{\circ} (t)
+ \di\frac{\partial \mathbb{L}}{\partial w} (\star_{\circ} (t)) \Delta_{\circ,+} \left [ \tilde{p}_{\circ} \right ] (t)
\right ) \,\Delta_{\frac{1}{2}}t ,
\end{array}
\right .
\end{equation}

where $\star_{\circ}=(p_\circ ,q_\circ ,\Delta_{\circ,+}[p_\circ],\Delta_{\circ,+}[q_\circ])$.\\

As the variation $\tilde{p}$ has no constraints, it remains a constant term when applying Discrete  integration by part formula and lemma \ref{dubois1}. Precisely, we obtain:

\begin{equation}
\left .
\begin{array}{lll}
D\mathscr{L}_{H,\T} (p, q) (\tilde{p},\tilde{q})
& = & \di\int_a^b 
\left ( 
\left [ 
\di\frac{\partial \mathbb{L}}{\partial q} (\star_{\circ}) 
\right ]_{\frac{1}{2},-}
- \Delta_{1/2,-} \left [ 
\di\frac{\partial \mathbb{L}}{\partial v} (\star_{\circ} ) \right ] 
\right )  (\sigma_{\circ}(t))
\tilde{q}(t) \,
\Delta t \\
& & + \di\int_a^b 
\left ( 
\left [
\di\frac{\partial \mathbb{L}}{\partial p} (\star_{\circ})
\right ]_{\frac{1}{2},-}
- \Delta_{1/2,-} 
\left [ 
\di\frac{\partial \mathbb{L}}{\partial w} (\star_{\circ} ) 
\right ] 
\right )  (\sigma_{\circ}(t))
\tilde{p} (t)
\, \Delta t \\
 & & + 
 \left [ 
\di\frac{\partial \mathbb{L}}{\partial w} (\star_{\circ} (t_{\frac{1}{2}})) \tilde{p}  (t_N)
-
\di\frac{\partial \mathbb{L}}{\partial w} (\star_{\circ} (t_\frac{1}{2})) \tilde{p} (t_0) 
 \right ]+ 
 \left [ 
\di\frac{\partial \mathbb{L}}{\partial p} (\star_{\circ} (t_{\frac{1}{2}})) \tilde{p} (t_N)
+
\di\frac{\partial \mathbb{L}}{\partial p} (\star_{\circ} (t_\frac{1}{2})) \tilde{p} (t_0) 
 \right ].
\end{array}
\right .
\end{equation}

Application of the discrete Dubois-Raymond lemma gives:
\begin{equation}
\left .
\begin{array}{l}
\left [\di \frac{\partial \mathbb{L}}{\partial p}(\star_{\circ} (t))
\right ]_{\frac{1}{2},-} -\Delta_{\frac{1}{2},-}\bigg[\di\frac{\partial \mathbb{L}}{\partial w}\bigg](\star_{\circ} (t)) 
=0, \\
\left [
\di\frac{\partial \mathbb{L}}{\partial q}(\star _\circ(t)) 
\right ]_{\frac{1}{2},-}
-\Delta_{\frac{1}{2},-}\bigg[\di\frac{\partial \mathbb{L}}{\partial v}\bigg] (\star_{\circ}(t)) 
=0 ,  
\end{array}
\right .
\end{equation}
for all $t\in \T_{\frac{1}{2}}^{\pm}$ and due to the fact that $\tilde{p}$ is free, we have also 
\begin{equation}
\left .
\begin{array}{lll}
\di\frac{\partial \mathbb{L}}{\partial p} (\star_{\circ} (t_{N-\frac{1}{2}})) 
+
\di\frac{\partial \mathbb{L}}{\partial w} (\star_{\circ} (t_{N-\frac{1}{2}}) )=0 ,\\
\di\frac{\partial \mathbb{L}}{\partial p} (\star_{\circ} (t_\frac{1}{2})) 
-
\di\frac{\partial \mathbb{L}}{\partial w} (\star_{\circ} (t_{\frac{1}{2}}))  
=0.
\end{array}
\right .
\end{equation}

As we have 
\begin{equation}
\frac{\partial \mathbb{L}}{\partial p}=v-\di\frac{\partial H}{\partial p},\ 
\frac{\partial \mathbb{L}}{\partial q}=-\di\frac{\partial H}{\partial q},\ 
\frac{\partial \mathbb{L}}{\partial w}=0,\  
\frac{\partial \mathbb{L}}{\partial v} =p ,
\end{equation}
we obtain

\begin{equation}
\label{main2}
\left .
\begin{array}{l}
\left [ \Delta_{\circ,+} [q_{\circ} ]-\di\frac{\partial H}{\partial p} (p_{\circ},q_{\circ} )
\right ]_{\frac{1}{2},-} 
=0, \\
\left [
-\frac{\partial H}{\partial q}(p_{\circ},q_{\circ}) 
\right ]_{\frac{1}{2},-}
-\Delta_{\frac{1}{2},-} [p_{\circ}] 
=0 .  
\end{array}
\right .
\end{equation}
for all $t\in \T_\frac{1}{2}^{\pm}$ and 
\begin{equation}
\label{bords}
\Delta_{\circ ,+}[q_{\circ}] (t_{N-\frac{1}{2}}) -\di\frac{\partial H}{\partial p} (q_{\circ} (t_{N-\frac{1}{2}} ),p_{\circ} (t_{N-\frac{1}{2}})) =0,\ \ 
\Delta_{\circ ,+}[q_{\circ}] (t_{\frac{1}{2}}) -\di\frac{\partial H}{\partial p} (q_{\circ} (t_{\frac{1}{2}} ),p_{\circ} (t_{\frac{1}{2}})) =0 .
\end{equation}

Equations \eqref{bords} can be used to simplify inductively the first equation of \eqref{main2}. Indeed, it follows that for all $t\in \T_{\frac{1}{2}}$, we have 
\begin{equation}
\Delta_{\circ,+} [q_{\circ} ]-\di\frac{\partial H}{\partial p} (p_{\circ},q_{\circ} ) =0 .
\end{equation}
\end{appendix}

\end{document}